\documentclass[10pt]{amsart}
\usepackage{amsmath}
\usepackage{latexsym}
\usepackage{amssymb}
\usepackage{amscd}
\input epsf

\usepackage{amsfonts}
\usepackage[all]{xy}

\DeclareFontShape{OT1}{cmr}{b}{n}{<12> cmr12}{}

\newtheorem{proposition}{Proposition}[section]

\newtheorem{lemma}[proposition]{Lemma}
\newtheorem{definition}[proposition]{Definition}
\newtheorem{theorem}[proposition]{Theorem}

\newtheorem{corollary}[proposition]{Corollary}

\newtheorem{remark}[proposition]{Remark}

\newenvironment{dok}{\par\vspace{-5pt}%
\par\noindent\begingroup%
\leftskip=0em\hspace{0em}{\bf Proof.}}%
{\endgroup\hfill$\Box$}

\newcounter{SNO}
\newcounter{THNO}[section]
\newcounter{tmp}

\def\db#1{ \bD^b({#1})}

\def\perf#1{{\mathfrak P}{\mathfrak e}{\mathfrak r}{\mathfrak f}({#1})}
\def\dsing#1{ \bD_{Sg}({#1})}
\def\dsingq#1{ \bD^{'}_{Sg}({#1})}
\def\lfr#1{{\mathfrak L}{\mathfrak f}{\mathfrak r}({#1})}
\def\dg{DG}
\def\dt{DB}

\def\abe{\operatorname{Pair}}
\def\ove{\overline}
\def\comp{\circ}
\def\Hom{{\mathrm H}{\mathrm o}{\mathrm m}}
\def\Hhom{{\mathbb H}{\mathrm o}{\mathrm m}}
\def\Cok{\operatorname{Cok}}

\def\h#1,#2{{\operatorname{Hom}}({#1}\:,\; {#2})}
\def\Ho#1,#2,#3,#4{{\operatorname{Hom}}^{#1}_{#2}({#3}\:,\; {#4})}
\def\Ex#1,#2,#3,#4{{\operatorname{Ext}}^{#1}_{#2}({#3}\:,\; {#4})}

\def\coh{\operatorname{coh}}
\def\Qcoh{\operatorname{Qcoh}}

\def\lto{\longrightarrow}

\def\wt{\widetilde}

\def\C{{\mathcal C}}

\def\D{{\mathcal D}}
\def\F{{\mathcal F}}
\def\E{{\mathcal E}}
\def\G{{\mathcal G}}

\def\N{{\mathcal N}}
\def\O{{\mathcal O}}
\def\R{{\mathcal R}}
\def\L{{\mathcal L}}

\def\M{{\mathcal M}}

\def\P{{\mathcal{P}}}

\def\ZZ{{\mathbb Z}}
\def\CC{{\mathbb C}}

\def\bD{{\mathbf D}}
\def\bR{{\mathbf R}}

\def\bL{{\mathbf L}}

\def\AA{{\mathbb A}}

\def\ZZ{{\mathbb Z}}

\def\CC{{\mathbb C}}

\def\PP{{\mathbb P}}

\def\Hom{\operatorname{Hom}}
\def\End{\operatorname{End}}
\def\Ext{\operatorname{Ext}}

\def\ext{\underline{{\mathcal E}xt}}
\def\rhom{\bR\underline{{\mathcal H}om}}
\def\hom{\underline{{\mathcal H}om}}

\def\Sing{\operatorname{Sing}}

\def\Ker{\operatorname{Ker}\,}
\def\Coker{\operatorname{Coker}\,}
\def\Im{\operatorname{Im}\,}

\def\Spec{\operatorname{Spec}}

\def\rH{{\operatorname{H}}}

\def\rZ{{\operatorname{Z}}}

\def\id{{\operatorname{id}}}

\def\mon{\rightarrowtail}
\def\epi{\twoheadrightarrow}

\def\depth{\operatorname{depth}}

\def\llto{{\;\relbar\joinrel\relbar\joinrel\relbar\joinrel\relbar\joinrel
\rightarrow\;}}

\title[]{Triangulated categories of singularities and
D-branes in Landau-Ginzburg models}

\author[]{Dmitri Orlov}

\address{ Algebra Section, Steklov Mathematical Institute RAN,
Gubkin str. 8, GSP-1, Moscow 119991, RUSSIA}

\email{orlov@mi.ras.ru}

\thanks{This work is done
under partial financial support of the Russian Foundation for
Basic Research (grant no.~02-01-00468), grant of President of RF
in support of young russian scientists No MD-2731.2004.1, grant
INTAS-OPEN-2000-269. The research describing in this work was made
possible in part by CRDF Award No RM1-2405-MO-02. It is also a
pleasure for me to express my gratitude to the Russian Science
Support Foundation.}

\date{}
\dedicatory{
Dedicated to the blessed memory of Andrei Nikolaevich Tyurin --
adviser and friend
}

\begin{document}

\maketitle

\vspace*{-0.2cm}
\section*{Introduction}

In spite of physics terms  in the title, this paper is purely
mathematical. Its purpose is to introduce triangulated categories
related to  singularities of algebraic varieties and establish a
connection of these categories with D-branes in Landau-Ginzburg
models
It seems that two different types of
categories can be associated with singularities (or singularities
of maps). Categories of the first type are connected with
vanishing cycles and closely related to the categories which were
introduced in \cite{Se} for symplectic Picard-Lefschetz pencils.
Categories of the second type are purely algebraic and come from
derived categories of coherent sheaves. Categories of this type
will be central in this work. An important notion here is the
concept of a perfect complex, which was introduced in \cite{Il}. A
perfect complex is a complex of sheaves which locally is
quasi-isomorphic to a bounded complex of locally free sheaves of finite
type (a good reference is \cite{TT}).

To any algebraic variety $X$ one can attach  the bounded derived
category of coherent sheaves $\db{\coh(X)}$. This category admits
a triangulated structure. The derived category of coherent sheaves
has a triangulated subcategory $\perf{X}$ formed by  perfect
complexes. If the variety $X$ is smooth then any coherent sheaf
has a finite resolution of locally free sheaves of finite type and
the subcategory of perfect complexes coincides with the whole of
$\db{\coh(X)}.$ But for singular varieties this property is not
fulfilled.
 We introduce a notion of triangulated category
of singularities $\dsing{X}$ as the quotient of the triangulated
category $\db{\coh(X)}$ by the full triangulated subcategory of
perfect complexes $\perf{X}$. The category $\dsing{X}$ reflects
the properties of the singularities of $X$ and "does not depend on
all of $X$". For example we prove that it is invariant with
respect to a localization in  Zariski topology (Proposition
\ref{locpr}). The category $\dsing{X}$ has  good properties when
$X$ is Gorenstein. In this case, if the locus of singularities is
complete then  all Hom's between  objects are finite-dimensional
vector spaces (Corollary \ref{fdho}).

The investigation  of  such  categories  is inspired by the
Homological Mirror Symmetry Conjecture (\cite{K}).

Works on topological string theory are mainly concerned  with the
case of  N=2 superconformal sigma-models with a
Calabi-Yau target space. In this case  the field theory has two
topologically twisted versions: A-models and B-models. The
corresponding D-branes are called A-branes and B-branes. The
mirror symmetry should interchange these two classes of D-branes.
From the mathematical point of view the category of B-branes on a
Calabi-Yau is the derived category of coherent sheaves on it
(\cite{K},\cite{Do}). As a candidate for a category of A-branes on
Calabi-Yau manifolds so-called Fukaya category has been proposed.
Its objects are, roughly speaking, Lagrangian submanifolds
equipped with flat vector bundles (\cite{K}). The Homological
Mirror Symmetry Conjecture asserts that if two Calabi-Yaus $X$ and
$Y$ are mirror to each other, then the derived category of
coherent sheaves on $X$ is equivalent to the Fukaya category of
$Y$, and vice versa.

On the other hand, physicists also consider more general
N=2  field theories and corresponding
D-branes. One class of such theories is given by a sigma-model
with a Fano variety as a target space. Another set of examples is
provided by N=2 Landau-Ginzburg models. In many cases these two
classes of N=2 theories are related by mirror symmetry
(\cite{HV}). For example, the sigma model with target $\PP^n$ is
mirror to a Landau-Ginzburg model which is given by superpotential
$$
W=x_1+\cdots +x_n + \frac{1}{x_1\cdots x_n}
$$
on $\CC^{*n}$. General definition of a Landau-Ginzburg model
involves, besides a choice of a target space, a choice of a
holomorphic function $W$ on this manifold. In particular,
non-trivial Landau-Ginzburg models require a non-compact target
space.

For Fano varieties one has the derived categories of coherent
sheaves (B-branes) and given  a symplectic form, one
can propose a suitable Fukaya category (A-branes).
Thus, if one wants to extend the Homological Mirror  Symmetry
Conjecture to the non-Calabi-Yau case, one should  understand
D-branes in Landau-Ginzburg models.

Categories of A-branes in Landau-Ginzburg models are studied in
\cite{HIV} and in \cite{Se} from the mathematical point of view.
Mirror symmetry relates B-branes on a Fano variety (coherent
sheaves) to A-branes in a LG model. In the case when the Fano
variety is $\PP^n$, the mirror symmetry prediction has been tested
in \cite{HIV} (and for $\PP^2$ in \cite{Se} from the mathematical
viewpoint).

One can also consider the  Fukaya category (A-branes) on a Fano
variety. One can expect that in this case the Fukaya category is
equivalent to the category of B-branes in the mirror
Landau-Ginzburg model. If we accept  mirror symmetry, we can
predict an answer  for the Fukaya category on a Fano variety by
studying B-branes in the mirror LG model. As a rule it is easier
to understand B-branes, rather than A-branes.

A mathematical definition for the category of B-branes in
Landau-Ginzburg models was proposed by M.Kontsevich. Roughly, he
suggests that  superpotential $W$ deforms complexes of coherent
sheaves to "twisted" complexes, i.e the composition of
differentials  is no longer zero, but  is equal to multiplication
by $W$. This "twisting" also breaks $\ZZ$\!-grading down to
$\ZZ/2$\!-grading. The equivalence of this definition with the
physics notion of B-branes in LG models was verified in  paper
\cite{KL} in  the case of the usual quadratic superpotential
$W=x_1^2+\cdots+ x_n^2$ and physical arguments were given
supporting Kontsevich's proposal for a general superpotential.

We establish a connection between  categories of B-branes in
Landau-Ginzburg models and triangulated categories of
singularities. We consider  singular fibres of the map $W$ and
show that the triangulated categories of singularities of these
fibres are equivalent to the categories of B-branes (Theorem
\ref{main2}). In particular this gives us that, in spite of the
fact that the category of B-branes is defined using the total
space $X$, it depends only on the singular fibres of the
superpotential. This result can be used for precise calculations
of the categories of B-branes in Landau-Ginzburg models. It is
remarkable fact that this construction was  known in the theory of
singularities as matrix factorization it was introduced in
\cite{Eis} for study of maximal Cohen-Macaulay modules over local
rings .

Another result which  is  useful for calculations connects the
categories of B-branes in different dimensions (Theorem
\ref{main1}). (This fact is  known in the local theory of
singularities as Kn\"orrer periodicity for maximal Cohen-Macaulay
modules.) It says the following: Let $W: X\to\CC$ be a
superpotential. Consider the manifold $Y=X\times \CC^2$ and
another superpotential $W'= W+xy$ on $Y$, where $x,y$ are
coordinates on $\CC^2$. Then the category of B-branes in the
Landau-Ginzburg model on $Y$ is equivalent to the category of
B-branes in the Landau-Ginzburg model on $X$. Actually, we prove
that the triangulated category of singularities of the fiber $X_0$
over point $0$ is equivalent to the triangulated category of
singularities of the fiber $Y_0$ (Theorem \ref{main1}). Keeping in
mind that  the categories of singularities are equivalent to the
categories of B-branes (Theorem \ref{main2}) we obtain the
connection between B-branes mentioned above.

In the end we give a calculation of the category of B-branes in
the Landau-Ginzburg model with the superpotential $W=z_0^n+
z_1^2+\cdots+z_{2k}^2$ on $\CC^{2k+1}$. The singularity of this
superpotential  corresponds to the Dynkin diagram $A_{n-1}$. This
category has $n-1$ indecomposable objects. We describe morphisms
between them, the translation functor and the triangulated
structure on this category.

\section{Singularities and triangulated categories}

\subsection{Triangulated categories and localizations}
In this section we remind definitions of a triangulated category
and its localization which were introduced in
\cite{Ve} (see also \cite{GM},\cite{KS},\cite{Ke}).
Let $\D$  be  an additive category.
The structure of {\sf a triangulated category} on $\D$
is defined by giving of the following data:
\begin{list}{\alph{tmp})}%
{\usecounter{tmp}}
\item an additive autoequivalence
$[1]: \D\lto\D$ (it is called a translation functor),
\item a class of exact (or distinguished) triangles:
$$
X\stackrel{u}{\lto}Y\stackrel{v}{\lto}Z\stackrel{w}{\lto}X[1],
$$
\end{list}
which must satisfy the set of axioms T1--T4.
\begin{itemize}
\item[T1.]
    \begin{itemize}
    \item[a)] For each object $X$ the triangle
    $X\stackrel{\id}{\lto}X\lto 0\lto X[1]$
 is exact.
    \item[b)] Each triangle isomorphic to an exact triangle
is exact.
    \item[c)] Any morphism $X\stackrel{u}{\lto} Y$ can be included
in an exact triangle\\
$X\stackrel{u}{\lto} Y\stackrel{v}{\lto} Z\stackrel{w}{\lto}
X[1]$.
    \end{itemize}
\item[T2.] A triangle
$X\stackrel{u}{\lto} Y\stackrel{v}{\lto} Z\stackrel{w}{\lto} X[1]$
is exact if and only if the triangle\\
$Y\stackrel{v}{\lto} Z\stackrel{w}{\lto}
X[1]\stackrel{-u[1]}{\lto} Y[1]$ is exact.
\item[T3.] For any two exact triangles and two morphisms
$f, g$ the diagram below
$$
\xymatrix{
X\ar[r]^u \ar[d]_f \ar@{}[dr]|
{\Box}&Y\ar[r]^v\ar[d]^g&Z\ar[r]^w\ar@{.>}[d]^h&X[1]
\ar[d]^{f[1]}\\
X'\ar[r]^{u'}&Y'\ar[r]^{v'}&Z'\ar[r]^{w'}&X'[1].
}
$$
can be completed to a morphism of triangles by a morphism $h:Z\to Z'$.
\item[T4.] For each pair of morphisms
$X\stackrel{u}{\lto} Y\stackrel{v}{\lto} Z$ there is a commutative
diagram
$$
\begin{CD}
X  @>{u}>>  Y @>{x}>> Z' @>>> X[1]\\
@| @V{v}VV   @VV{w}V  @| \\
X   @>>>    Z  @>{y}>>     Y'   @>{w'}>>  X[1]\\
@.   @VVV       @VV{t}V  @VV{u[1]}V\\
@. X' @=  X'  @>{r}>>   Y[1]\\
@. @VV{r}V      @VVV   @.      \\
@. Y[1]  @>{x[1]}>>    Z'[1] @.
\end{CD}
$$
where the first two rows and the two central columns are exact triangles.
\end{itemize}

A functor $F : {\D} \lto{\D}'$ between two triangulated categories
${\D}$ and ${\D}'$ is called {\sf exact} if it commutes with the
translation functors, i.e there fixed a natural isomorphism
$
t : F\circ [1]_{\D} \stackrel{\sim}{\to}[1]_{\D'} \circ F
$
and
it  transforms all exact triangles into exact
triangles, i.e for each exact triangle
 $X\stackrel{u}{\to} Y\stackrel{v}{\to} Z\stackrel{w}{\to}X[1]$
in ${\D}$
the triangle
$$
FX\stackrel{Fu}{\lto} FY\stackrel{Fv}{\lto} FZ\stackrel{t_X Fw}{\lto}
FX[1]
$$
is exact.

A full additive subcategory  $\N\subset \D$ is called full
triangulated subcategory, if the following conditions hold: it is
closed with respect to the translation functor in $\D$ and if it
contains any two objects of an exact triangle in $\D$ then it
contains the third object of this triangle as well. The full
triangulated subcategory $\N\subset \D$ is called {\sf thick} if
it is closed with respect to taking of direct summands in $\D$.

Now we remind the definition of a localization of categories. Let
$\C$ be a category and let $\Sigma$ be a class of morphisms in
$\C$. It is well-known \cite{GZ} that there is a large category
$\C[\Sigma^{-1}]$ and a functor $Q:\C\to \C[\Sigma^{-1}]$ which is
universal among the functors making the elements of $\Sigma$
invertible. (Note that the objects of the category $\C[\Sigma^{-1}]$ form not a set, but a class in general.)
The category $\C[\Sigma^{-1}]$ has a good description
if $\Sigma$ is a multiplicative system.

A family of morphisms $\Sigma$ in a category $\C$
is called {\sf a multiplicative system} if it
satisfies the following conditions:
\begin{itemize}
\item[M1.] all identical morphisms $\id_{X}$
belongs to
$\Sigma$;
\item[M2.] the composition of two elements of
$\Sigma$ belong to
$\Sigma$;
\item[M3.] any diagram
$X'\stackrel{s}{\longleftarrow} X\stackrel{u}{\lto} Y$,
with $s\in\Sigma$ can be completed to a commutative square
$$
\xymatrix{
X \ar[r]^{u}\ar[d]_{s} & Y \ar@{-->}[d]^{t}\\
X' \ar@{-->}[r]^{u'}& Y'}
$$
with $t\in \Sigma$ (the same when all arrows reversed);
\item[M4.] for any two morphisms $f,g$ the existence
of $s\in\Sigma$ with
$fs=gs$ is equivalent to the existence of $t\in \Sigma$
with $tf=tg$.
\end{itemize}

If $\Sigma$ is a multiplicative system then  $\C[\Sigma^{-1}]$
has the following description.
The objects of $\C[\Sigma^{-1}]$ are the objects of
$\C$.
The morphisms  from $X$  to $Y$ in
$\C[\Sigma^{-1}]$ are  pairs
$(s,f)$ in $\C$
of the form
 $$
X\stackrel{f}{\lto}Y'\stackrel{s}{\longleftarrow}  Y,
\qquad s\in \Sigma
$$
 modulo the following equivalence relation:
 $(f,s)\sim (g,t)$
iff there is a commutative diagram
$$
\xymatrix{
&Y'\ar[d]&\\
X\ar[ur]^f \ar[r]^{h} \ar[dr]_{g}& Y''' &Y\ar[ul]_s \ar[l]_{r}\ar[dl]^t\\
&Y''\ar[u]}
$$
with $r\in\Sigma$.

The composition of the morphisms
$(f,s)$  and $(g,t)$ is a morphism $(g'f, s't)$
defined from the following diagram, which exists by M3:
$$
\xymatrix{
&&Z''&&\\
&Y'\ar@{-->}[ur]^{g'}&&Z'\ar@{-->}[ul]_{s'}&\\
X\ar[ur]^{f}&&Y\ar[ul]_{s}\ar[ur]^g && Z\ar[ul]_t}.
$$

It can be checked that $\C[\Sigma^{-1}]$ is a category and there
is a quotient functor
$$
Q: \C \lto \C[\Sigma^{-1}], \quad X\mapsto X, f \mapsto (f,1)
$$
which inverts all elements of $\Sigma$
and it is universal in this sense (see \cite{GZ}).

Let $\D$ be a triangulated category and $\N\subset \D$
be a full triangulated subcategory.
Denote by $\Sigma(\N)$ a class of morphisms $s$ in $\D$
embedding into an exact triangle
$$
X\stackrel{s}{\lto} Y\lto N\lto X[1]
$$
with $N\in \N$.
It can be checked that $\Sigma(N)$ is a multiplicative system.
We define
$$
\D/\N:=\D[\Sigma(\N)^{-1}].
$$
We endow the category $\D/\N$ with a translation functor
induced by the translation functor in the category $\D$.
\begin{lemma}
The category $\D/\N$ becomes a triangulated category by taking for
exact triangles such that are isomorphic to the images of exact
triangles in $\D$. The quotient functor $Q:\D\lto \D/\N$
annihilates $\N$. Moreover, any exact functor $F: \D\lto \D'$ of
triangulated categories for which $F(X)\simeq 0$ when $X\in \N$
factors uniquely through $Q$.
\end{lemma}
The following lemma, which will be necessary in the future, is evident.
\begin{lemma}\label{adjqu}
Let $\M$ and $\N$ be full triangulated subcategories
in  triangulated categories $\C$ and $\D$ respectively.
Let $F: \C\to \D$ and $G: \D\to \C$ be  adjoint pair
of exact functors such that $F(\M)\subset \N$ and $G(\N)\subset \M$.
Then they induce functors
$$
\ove{F}:\C/\M\lto\D/\N,
\qquad
\ove{G}:\D/\N\lto \C/\M
$$
which are adjoint.
\end{lemma}

\begin{proposition}{\rm (\cite{Ve}, \cite{KS}).}
\label{ffemb}
Let $\D$ be a triangulated category and $\D', \N$ be full
triangulated subcategories. Let $\N'=\D'\cap \N$.
Assume that any morphism $N\to X'$ with $N\in \N$ and $X'\in D'$
admits a factorization $N\to N'\to X'$ with $N'\in \N'$.
Then the natural functor
$$
\D'/\N'\lto \D/N
$$
is fully faithful.
\end{proposition}

\subsection{Triangulated categories of singularities}

Let $X$ be a scheme
over a field $k$. We will say that it satisfies condition (ELF) if it is

\vspace{0pt}
\begin{tabular}{ll}
(ELF)&\quad \begin{tabular}{l}
 separated,  noetherian, of finite  Krull dimension, and  the category of coherent sheaves\\
$\coh(X)$ has enough  locally free sheaves.\\
\end{tabular}
\end{tabular}
\vspace{3pt}

The last condition means that for any coherent sheaf $\F$
there is a vector bundle $\E$ and an epimorphism $\E\twoheadrightarrow\F.$
For example any quasi-projective scheme satisfies these conditions.
Note that any closed and any open subscheme of $X$ is also
noetherian, finite dimension and has enough locally free sheaves.
It is clear for a closed subscheme while for an open subscheme
$U$ it follows from the fact that any coherent sheaf on $U$
can be obtained as the restriction of a coherent sheaf on $X$ (see
\cite{H}, ex.5.15).

Denote by $\db{\coh(X)}$ (resp. $\db{\Qcoh(X)}$)
 the bounded derived categories of coherent (resp. quasi-coherent)
 sheaves on $X$.
 These  categories  have  canonical triangulated structures.

Since $X$ is noetherian the natural functor $\db{\coh(X)}\lto
\db{\Qcoh(X)}$ is fully faithful and realizes an equivalence of
$\db{\coh(X)}$ with the full subcategory
$\db{\Qcoh(X)}_{\coh}\subset \db{\Qcoh(X)}$ consisting of all
complexes with coherent cohomologies (see \cite{Il} II, 2.2.2). By
this reason,  when we  consider $\db{\coh(X)}$ as a subcategory of
$\db{\Qcoh(X)}$ we will identify it with the full subcategory
$\db{\Qcoh(X)}_{\coh}$, adding all isomorphic objects.

\begin{lemma}{\rm(\cite{Il}, \cite{TT})}\label{cov}
Let $X$ be as above. Then for any bounded above complex  of
quasi-coherent sheaves $C^{\cdot}$ on $X$ there is a bounded above
complex of locally free sheaves  $P^{\cdot}$ and a
quasi-isomorphism of the complexes $P^{\cdot}\stackrel{\sim}{\lto}
C^{\cdot}$. Moreover, if $C^{\cdot}\in \db{\Qcoh(X)}_{\coh}$ then
there is a bounded above complex of locally free sheaves of finite
type $P^{\cdot}$ with a quasi-isomorphism
$P^{\cdot}\stackrel{\sim}{\lto}C^{\cdot}$.
\end{lemma}

Recall the constructions of standard truncation functors. Let
$C^{\cdot}$ be a complex. There is brutal truncation
$$
\sigma^{\ge k}C^{\cdot}=\cdots \to 0\to C^k\to C^{k+1}\to C^{k+2}\to \cdots
$$
This is a subcomplex of $C^{\cdot}$. The quotient
$C^{\cdot}/\sigma^{\ge k}C^{\cdot}$ is another brutal truncation
$\sigma^{\le k-1}C^{\cdot}$.

There is also the good truncation
$$
\tau^{\ge k}C^{\cdot}=\cdots\to 0\to \Im d\; C^{k-1}\to C^{k}\to
C^{\cdot}\to\cdots
$$
There is a quotient map $C^{\cdot}\twoheadrightarrow \tau^{\ge
k}C^{\cdot}$ which induces  isomorphisms on cohomologies $H^n$ for
all $n\ge k$. The kernel of this map is denoted $\tau^{\le
k-1}C^{\cdot}$.

\begin{definition}
A bounded complex of coherent
sheaves will be called {\sf a perfect complex} if it is quasi-isomorphic
to a bounded complex
of locally free sheaves of finite type.
\end{definition}
\begin{lemma}\label{locper}
Any complex $C^{\cdot},$ which is isomorphic to a bounded complex
of locally free sheaves in $\db{\coh(X)},$ is perfect.
\end{lemma}
\begin{dok}
We can represent the isomorphism in the derived category via
calculus of fractions as $P^{\cdot}\stackrel{s}{\leftarrow}
E^{\cdot} \stackrel{t}{\to} C^{\cdot}$, where $P^{\cdot}$ is a
bounded complex of locally free sheaves and $s, t$ are
quasi-isomorphisms. By Lemma \ref{cov}, there is a bounded above
complex $Q^{\cdot}$ of locally free sheaves and quasi-isomorphism
$ Q^{\cdot}\to E^{\cdot}$. Consider a good truncation $\tau^{\ge
-k}Q^{\cdot}$ for sufficiently large $k$. As $E^{\cdot}$ is
bounded there is a morphism $ r:\tau^{\ge -k}Q^{\cdot}\to
E^{\cdot}$ that is also a quasi-isomorphism. To prove the lemma it
is sufficient to show that $\tau^{\ge -k}Q^{\cdot}$ is a complex
of locally free sheaves. Consider the composition $sr:\tau^{\ge
-k}Q^{\cdot}\to P^{\cdot}$ which is a quasi-isomorphism. The cone
of $sr$ is a bounded acyclic complex all terms of which, excepted
maybe  the leftmost term, are locally free. This implies that the
leftmost term is locally free as well, because the kernel of an
epimorphism of locally free sheaves is locally free. Thus
$\tau^{\ge -k}Q^{\cdot}$ is a bounded complex of locally free
sheaves and hence $C^{\cdot}$ is perfect.
\end{dok}

The perfect complexes form a full triangulated subcategory
$\perf{X}\subset \db{\coh(X)}$, which is thick.
\begin{remark}
Actually, a perfect complex is defined as a complex of
$\O_X$\!-modules locally quasi-isomorphic to a bounded complex of
locally free sheaves of finite type. But under our assumption on
the scheme any such complex is quasi-isomorphic to a bounded
complex of locally free sheaves of finite type(see \cite{Il} II,
or \cite{TT} \S 2).
\end{remark}
\begin{definition} We define a triangulated category
$\dsing{X}$ as the quotient
of the triangulated category $\db{\coh(X)}$ by the full triangulated
subcategory
$\perf{X}$ and call it as a triangulated category of singularities of $X$.
\end{definition}

\begin{remark}
{\rm It is  known that if a scheme $X$ is as above and is
regular in addition then the
subcategory of perfect complexes coincides with the whole bounded
derived category of coherent sheaves. Hence, the triangulated category
of singularities is trivial in this case.}
\end{remark}

We also consider  the full
triangulated subcategory $\lfr{X}\subset \db{\Qcoh(X)}$ consisting of objects which are
isomorphic to bounded complexes of locally free sheaves in
$\db{\Qcoh(X)}$. By the same argument as in Lemma \ref{locper} we
can show that for any  complex $C^{\cdot}\in \lfr{X}$ there is a
bounded complex of locally free sheaves $P^{\cdot}$ and a
quasi-isomorphism $P^{\cdot}\to C^{\cdot}$. Using Lemma \ref{cov},
it can be  checked that the subcategory $\db{\coh(X)}\cap\lfr{X}$
coincides with the subcategory of perfect complexes $\perf{X}$.

Define a triangulated
category $\dsingq{X}$ as
 the quotient $\db{\Qcoh(X)}/\lfr{X}$.
The full embedding $\db{\coh(X)}\lto \db{\Qcoh(X)}$
induces a  functor $\dsing{X}\lto \dsingq{X}$.
We will show that this functor is also full embedding.

\begin{lemma}\label{persup}
 Let $X$ be as above and let $\F$ be a
quasi-coherent sheaf on $X$ such that for any point $x\in
\Sing(X)$ it is locally free in some neighborhood of $x$. Then it
belongs to the subcategory $\lfr{X}$. If $\F$ is coherent in
addition then it is perfect as a complex.
\end{lemma}
\begin{dok}
By Lemma \ref{cov}, there is a bounded above complex $Q^{\cdot}$
of locally free sheaves and quasi-isomorphism $ Q^{\cdot}\to
\F$.Consider a good truncation $\tau^{\ge -k}Q^{\cdot}$ for
sufficiently large $k$. There is a morphism $ r:\tau^{\ge
-k}Q^{\cdot}\to \F$ that is also a quasi-isomorphism. To prove the
lemma it is sufficient to show that $\tau^{\ge -k}Q^{\cdot}$ is a
complex of locally free sheaves. All terms of this complex, except
maybe for the leftmost term, are locally free. But for any point
$x\in \Sing(X)$ the leftmost term is also locally free in some neighborhood of $x$, because
$\F$ is locally free there. If now $x\not\in\Sing(X)$ then there is a
neighborhood $U$ of $x$ which is regular. Hence, the leftmost term
is locally free on $U$ under assumption that $k>\dim X$.

If now $\F$ is coherent then we can take the resolution $Q^{\cdot}$
such that all terms are of finite type. Hence, $\F$ is perfect.
\end{dok}

In particular, we get from this lemma  that if the support
of an object $E^{\cdot}\in \db{\coh(X)}$ does not meet $\Sing(X)$
then $E^{\cdot}$ is  perfect.

\begin{lemma}\label{lfint}
Let $X$ be as above.
Then any object $A\in\dsing{X}$ is isomorphic
to an object $\F[k]$
where $\F$ is a coherent sheaf.
\end{lemma}
\begin{dok} The object $A$ is a bounded complex of coherent sheaves.
Let us take locally free  bounded above resolution
$P^{\cdot}\stackrel{\sim}{\to} A$
which exists by Lemma \ref{cov}.
Consider a brutal truncation $\sigma^{\ge -k}P^{\cdot}$
for sufficiently large $k\gg 0$. Denote by $\F$ the cohomology
$H^{-k}(\sigma^{\ge -k}P^{\cdot})$.
It is clear that $A\cong \F[k+1]$ in $\dsing{X}$.
\end{dok}

\begin{lemma}\label{van}
Let $X$ be as above. Then for any locally free sheaf
$\E$ and for any quasi-coherent sheaf $\F$
$$
\Ext^i(\E, \F)=0,\qquad
\text{for}
\quad i>n.
$$
\end{lemma}
\begin{dok}
Let $U_1 \cup ... \cup U_n$ be an affine cover of $X$.
For all subsequence of indices $I=(i_1,..., i_k)$,
let $U_{I}=U_{i_1}\cap ... \cap U_{i_k},$ and
let $j_{I}: U_{I}\hookrightarrow X$ be the open immersion.
As $X$ is separated, each $U_{I}$ is affine and each $j_{I}$
is an affine map. Hence, $j_{I*}$ is exact and preserves
quasi-coherence.

We consider the \u{C}ech hypercover complex of quasi-coherent sheaves
\begin{equation}\label{cmpl}
0\to \F\to\bigoplus_{i=1}^{n} j_{i*}j^*_{i} \F\to
\bigoplus_{I=(i_1, i_2)} j_{I*}j^*_{I} \F\to
\bigoplus_{I=(i_1, i_2, i_3)} j_{I*}j^*_{I} \F\to\cdots
\end{equation}
This is an exact sequence of sheaves.

Since $U_I=\Spec(A_I)$ is affine the category
of quasi-coherent sheaves on $U_I$ is equivalent to the category
of $A_I$\!-modules. Therefore, we have
 $$
\Ext^i(\E, j_{I*}j^*_{I}\F)=
\Ext^i(j_{I}^*\E, j^*_{I}\F)=0
\qquad\text{for all}
\qquad i>0,
$$
because $j_{I}^*\E$ corresponds to a projective module over $A_I$.
Thus the nontrivial Ext's from $\E$ to $\F$ are bounded by the length
of the complex (\ref{cmpl}) which is equal to $n$.
\end{dok}

\begin{proposition}\label{femb}
For a scheme $X$ as above the natural functor
$\dsing{X}\lto\dsingq{X}$ is fully faithful.
\end{proposition}
\begin{dok}
Consider the full embedding $\db{\coh(X)}\hookrightarrow
\db{\Qcoh(X)}$. We have $\perf{X}=\db{\coh(X)}\cap \lfr{X}$. To
prove the proposition we check that the conditions of Proposition
\ref{ffemb} are fulfilled. Let $Q^{\cdot}\stackrel{t}{\to}
\F^{\cdot}$ be a morphism in $\db{\Qcoh(X)}$ such that
$Q^{\cdot}\in \lfr{X}$ and $\F^{\cdot}\in \db{\coh(X)}$. Let
$P^{\cdot}\stackrel{s}{\to}\F^{\cdot}$ be a quasi-isomorphism
where $P^{\cdot}$ is a bounded above complex of locally free
sheaves of finite type. It exists by Lemma \ref{cov}. The brutal
truncation $\sigma^{\ge -k}P^{\cdot}$ is a perfect complex.
Consider the map $r:\sigma^{\ge -k}P^{\cdot}\to \F^{\cdot}$
induced by $s$. A cone of the map $r$ is isomorphic to $G[k+1]$
where $G$ is a coherent sheaf. Since $Q^{\cdot}$ is a bounded
complex of locally free sheaves by Lemma \ref{van}, we obtain that
$$
\Hom(Q^{\cdot}, G[k+1])=0
$$
for sufficiently large $k$.
Hence, the map $t$ can be lifted to some map
$Q^{\cdot}\to \sigma^{\ge -k}P^{\cdot}$.
\end{dok}

Consider a morphism $f: X\to Y $
of finite Tor-dimension (for example flat). It  defines
an inverse image functor $\bL f^*:\db{\coh(Y)}\lto \db{\coh(X)}$.
It is clear that the functor $\bL f^*$ sends
a perfect complex on $Y$ to a perfect complex on $X$.
Therefore, we get a functor
$\bL \bar{f}^* :\dsing{Y} \lto \dsing{X}$.
By the same reason there is a functor $\bL \bar{f}^* : \dsingq{Y}
\lto\dsingq{X}$.

Let  $f: X\to Y$ be a morphism of finite Tor-dimension.
Suppose that it is  proper
of locally finite type. Then there is the functor
of direct image $\bR f_* : \db{\coh(X)}\to \db{\coh(Y)}$
and, moreover, this functor takes  perfect complexes to  perfect
complexes (see \cite{Il} III, or \cite{TT}). In this case we get a functor
$\bR \bar{f}_* :\dsing{X}\to \dsing{Y}$, which is the right
adjoint to $\bL \bar{f}^*$.

Now we prove a local property for triangulated categories of
singularities.
\begin{proposition}\label{locpr}
Let X be as above and let $j: U\hookrightarrow X$ be an embedding
of an open subscheme such that $\Sing(X)\subset U$. Then the
functor $\bar{j}^*:\dsing{X}\lto \dsing{U}$ is an equivalence of
triangulated categories.
\end{proposition}
\begin{dok}
Since $X$ is noetherian there is the functor $\bR j_* :
\db{\Qcoh(U)}\lto \db{\Qcoh(X)}$ which is right adjoint to $j^*$.
The composition $j^* \bR j_*$ is isomorphic to the identity
functor. Take an object $B\in \lfr{U}$ and consider $\bR j_*(B)$.
It is easy to see that the object $\bR j_*(B)$ belongs to
$\lfr{X}$. Actually, this condition is local. For $U$ it is
fulfilled and for $X\setminus \Sing(X)$ as for smooth scheme it is
evident. Thus the functor $\bR j_*$ induces the functor
$$
\bR \bar{j}_*: \dsingq{U}\lto \dsingq{X}.
$$
Moreover, the functor $\bR \bar{j}_*$ is right adjoint to $\bar{j}^*$.

For any object $A\in \db{\Qcoh(X)}$ we have a canonical
map
$
\mu_A: A\lto \bR j_* j^* A
$.
A cone $C(\mu_A)$ of this map is an object whose support
 belongs to $X\setminus U$ and does not intersect
$\Sing(X)$. Hence, by Lemma \ref{persup} the object $C(\mu_A)$
belongs to the subcategory $\lfr{X}$. This gives that $\mu_A$
becomes an isomorphism in $\dsingq{X}$. Therefore, the functor
$$
\bar{j}^*:\dsingq{X}\lto
\dsingq{U}
$$
is fully faithful.
On the other hand, we know that $j^* \bR j_*(B)\cong B$
for any $B\in \db{\Qcoh(U)}$. Hence, $\bar{j}^*$ is an equivalence.

The functor $j^*$ preserves coherence. Thus, using
Proposition \ref{femb}, we obtain that the functor
$$
\bar{j}^*:\dsing{X}\lto
\dsing{U}
$$
is fully faithful. Now note that by Lemma \ref{lfint} any object
$B\in \dsing{U}$
is isomorphic to $\F[k]$ where $\F$ is a coherent sheaf on $U$
and any coherent sheaf on $U$ can be obtained as the restriction of
a coherent sheaf on $X$ (\cite{H}, Ex.5.15).
This implies that $\bar{j}^*: \dsing{X}\lto
\dsing{U}$ is an equivalence.
\end{dok}

\subsection{Triangulated categories of singularities
for Gorenstein schemes}

Remind the definition of a Gorenstein local ring and
a Gorenstein scheme.
\begin{definition}
A local noetherian ring $A$ is called Gorenstein if
$A$ as a module over itself has a finite injective resolution.
\end{definition}
It can be shown that if $A$ is Gorenstein than $A$ is a
dualizing complex for itself (see \cite{Ha}).
This means that $A$ has  finite injective dimension and
the natural map
$$
M\lto \bR\Hom^{\cdot}(\bR\Hom^{\cdot}(M, A), A)
$$
is an isomorphism for any coherent $A$\!-module $M$ and as
consequence for any object from $\db{\coh(\Spec(A))}$.

\begin{definition} A scheme $X$ is  Gorenstein if
all of its local rings are Gorenstein local rings.
\end{definition}
\begin{remark}\label{finex}
{\rm If $X$ is Gorenstein and has finite Krull dimension, then
$\O_X$ is a dualizing complex for $X$, i.e. it has finite injective
dimension as quasi-coherent sheaf and the natural map
$$
\F\lto \rhom^{\cdot}(\rhom^{\cdot}(\F, \O_X), \O_X)
$$
is an isomorphism for any coherent sheaf $\F$.
In particular, we have that  there is
an integer $n_0$ such that
$\ext^{i}(\F, \O_X)=0$ for each quasi-coherent sheaf
$\F$ and all $i>n_0$.}
\end{remark}

\begin{lemma}\label{bext}
 Let $X$ be as above and Gorenstein.
Then for any coherent sheaf $\F$ and
an object $P^{\cdot}\in \perf{X}$ there is an integer $m$
depending only on $P^{\cdot}$ such that
$$
\Hom^{i}(\F, P^{\cdot})=0
$$
for all $i>m$.
\end{lemma}
\begin{dok}
For any coherent sheaf $\F$ and for any locally free sheaf $\P$
we know that
$$
\ext^{i}(\F, \P)=\ext^{i}(\F, \O_X)\otimes \P=0
$$
for any $i>n_0$.
Using the spectral sequence from local to global Ext's
we obtain that
$$
\Ext^{i}(\F, \P)=0
$$
for $i>n_0+n$ where $n$ is dimension of $X$.
Since  $P^{\cdot}$ is a bounded complex
there exists $m$ depending only on $P^{\cdot}$
such that
$$
\Hom^{i}(\F, P^{\cdot})=0
$$
for $i>m$.
\end{dok}

\begin{lemma}\label{rres}
Let X as above and Gorenstein. Then the following
conditions on a coherent sheaf $\F$ are equivalent.
\begin{itemize}
\item[1)] The sheaves $\ext^{i}(\F, \O_X)$ are trivial for all
$i>0$.

\item[2)] There is a right locally free resolution
$$
0\lto\F\lto\{Q^0\lto Q^1 \lto Q^3\cdots \}.
$$
\end{itemize}
\end{lemma}
\begin{dok}
$1)\Rightarrow 2).$
Denote by $\F^{\vee}$ the sheaf $\hom(\F, \O_X)$.
Consider a left locally free resolution of $\F^{\vee}$.
Applying to it the functor $\hom( \cdot, \O_X)$ we obtain a right
locally free resolution of $\F$, because $\O_X$ is dualizing complex.

$2)\Rightarrow 1).$
Consider a brutal truncation $\sigma^{\le k}Q^{\cdot}$
for sufficiently large $k$.
Denote by $\G$ the nontrivial cohomology $H^k(\sigma^{\le k}Q^{\cdot})$.
For any $i>0$ we have  isomorphisms
$$
\ext^{i}(\F, \O_X)\cong \ext^{i+k+1}(\G, \O_X)=0
$$
The last equality follows from Remark \ref{finex}.
\end{dok}

\begin{lemma}\label{locfr}
Let X be as above and Gorenstein.
Let $\F$ be a coherent sheaf which is perfect as a complex.
Suppose that
$
\ext^{i}(\F, \O_X)=0$
for all $i>0$. Then $\F$ is locally free.
\end{lemma}
\begin{dok}
If $\F$ is perfect then $\F^{\vee}=\hom(\F, \O_X)$
is also perfect. Hence, it has bounded locally free resolution
$\P^{\cdot}\lto\F^{\vee}$. Since $\O_X$ is a dualizing complex
we obtain a bounded right resolution $\F\lto \P^{\cdot\vee}$.
Thus $\F$ is locally free.
\end{dok}

\begin{proposition}\label{stab}
Let $X$ be as above and  Gorenstein. Let
$\F$ and $\G$ be  coherent sheaves such that
$\ext^{i}(\F, \O_X)=0$ for all $i>0$. Fix $N$ such that
$\Ext^{i}(\P, \G)=0$ for $i>N$ and for any locally free sheaf $\P$.
Then
$$
\Hom_{\dsing{X}}(\F, \G[N])\cong \Ext^{N}(\F, \G)/\R
$$
where $\R$ is the subspace of elements factoring through locally free,
i.e.  $e\in \R$ iff $e=\beta\alpha$ with $\alpha:\F\to \P$
and $\beta\in \Ext^{N}(\P, \G)$ where $\P$ is locally free.
\end{proposition}
\begin{dok}
By the definition of a localization any morphism from $\F$ to $\G[N]$
in $\dsing{X}$ can be represented by a pair of morphisms in $\db{\coh(X)}$
of the form
\begin{equation}\label{domik}
\F\stackrel{s}{\longleftarrow}A\stackrel{a}{\lto} \G[N]
\end{equation}
such that the cone $C(s)$ is a perfect complex.
By  Lemma \ref{rres} there is a right locally free resolution
$\F\lto Q^{\cdot}$. We consider a brutal truncation
$\sigma^{\le k} Q^{\cdot}$ for sufficiently large $k$
such that $\Hom(\E[-k-1], C(s))=0,$ where $\E=H^{k}(\sigma^{\le k}Q^{\cdot})$.
Such $k$ exists by Lemma \ref{bext}.
Using the triangle
$$
\E[-k-1]\lto \F\lto \sigma^{\le k}Q^{\cdot}\lto \E[-k]
$$
we find that the map $\F\to C(s)$ can be lifted to a map
$\sigma^{\le k}Q^{\cdot}\to C(s)$.
Therefore, there is a map $\E[-k-1]\to A$ which induces a pair of the form
\begin{equation}\label{domik2}
\F\stackrel{s'}{\longleftarrow}\E[-k-1]\stackrel{e}{\lto} \G[N],
\end{equation}
and this pair gives the same morphism in $\dsing{X}$ as the pair (\ref{domik}).
Since $\Ext^{i}(\P, \G)=0$ for $i>N$ and for any locally free sheaf
$\P$, we obtain
$$
\Hom(\sigma^{\le k}Q^{\cdot}[-1], \G[N])=0.
$$
Hence, there is a map $f$ which completes the following diagram
$$
\xymatrix{
&\E[-k-1]\ar[dr]^{e}\ar[dl]_{s'}\\
\F\ar[rr]^{f}&& \G[N]}.
$$
This gives us that the map $f$ is
equivalent
to the map (\ref{domik2}) and as consequence to
the map (\ref{domik}).
Thus, any morphism from $\F$ to $\G[N]$ in $\dsing{X}$
is represented by  a morphism from $\F$ to $\G[N]$ in
the category $\db{\coh(X)}$.
Now if $f$ is the $0$\!-morphism in $\dsing{X}$ then
repeating the procedure above we find that $f$ factors
through a morphism $\sigma^{\le k}Q^{\cdot}\to \G[N]$.
By the condition on $G$ any such morphism can be lifted
to a morphism $Q^0\to \G[N]$.
Hence, if $f$ is  the $0$\!-morphism in $\dsing{X}$
then it factors through $Q^0$.
\end{dok}

\begin{remark} {\rm If $X$ is affine
 then $N$ can be taken equal to $0$.}
\end{remark}

\begin{proposition}\label{fint}
Let $X$ be as above and Gorenstein.
Then any object $A\in\dsing{X}$ is isomorphic
to the image of a coherent sheaf $\F$ such that $\ext^{i}(\F, \O_X)=0$
for all $i>0$.
\end{proposition}
\begin{dok} An object $A$ is a bounded complex of coherent sheaves.
Let us take locally free  bounded above resolution
$P^{\cdot}\stackrel{\sim}{\to} A$
which exists by Lemma \ref{cov}.
Consider a brutal truncation $\sigma^{\ge -k}P^{\cdot}$
for sufficiently large $k\gg 0$. Denote by $\G$ the cohomology
$H^{-k}(\sigma^{\ge -k}P^{\cdot})$.
Since $A$ is bounded and $X$ is Gorenstein we know that the complex
$\rhom^{\cdot}(A, \O_X)$ is bounded.
This implies that if  $k\gg 0$ then
$\ext^{i}(\G, \O_X)=0$ for all $i>0$.
Moreover, we get that $A\cong \G[k+1]$ in $\dsing{X}$.

By Lemma \ref{rres} there is a right locally free resolution
$\G\to {Q^{0}}\to {Q^{1}}\to \cdots$.
Consider $\Im d_{k-1}=\Ker d_k\subset {Q^{k}}$ and denote
it by $\F$. Applying again Lemma \ref{rres} we obtain that
 $\ext^{i}(\F, \O_X)=0$ for
all $i\ge 0$.
And we have an isomorphism
$A\cong \G[k+1]\cong \F$
in $\dsing{X}$.
\end{dok}

\begin{corollary}\label{fdho}
Let $X$ be as above and Gorenstein such that  the closed subset $\Sing(X)$ is complete.
Then $\dim_{k}\Hom(A, B)<\infty$ for any two objects $A, B\in
\dsing{X}$.
\end{corollary}
\begin{dok}
Let $\overline{X}$ be  some compactification of $X$. Since
$\Sing(X)$ is complete the intersection $\Sing(X)$ with the
complement $\overline{X}\backslash X$ is empty. Resolving if it is
necessary singularities of $\overline{X}$ on the complement
$\overline{X}\backslash X$ we can assume that
$\Sing(\overline{X})$ coincides with $\Sing(X)$. By  Proposition
\ref{locpr} there is an equivalence
$\dsing{X}\simeq\dsing{\overline{X}}$.
 We know that for any two objects of
$\db{\coh(\overline{X})}$ the space of morphisms is finite
dimensional. Now the statement of the corollary immediately
follows from Propositions \ref{stab} and \ref{fint}.
\end{dok}

\section{Kn\"orrer periodicity}

Let $X$ be a separated regular noetherian scheme of finite
Krull dimension.
Any such scheme has enough locally free sheaves
(see \cite{Il}, II).
Let $f: X\to \AA^1$
be a flat morphism. Consider the scheme $Y=X\times \AA^2$
and a morphism $g=f+xy$ to $\AA^1$, where $x, y$ are coordinates
on $\AA^2$.
Denote by $X_0=X/f$ and $Y_0=Y/g$ the fibers of $f$ and $g$ respectively
over the point $0$. Consider the scheme $Z =Y_0/x$.
There are  natural maps $i: Z\to Y_0$ and $q: Z\to X_0$
where the former is a closed embedding and the latter
is an $\AA^1$\!-bundle.
All the schemes, introduced above, are separated noetherian
schemes of finite Krull dimension that have enough locally free
sheaves.

Consider the composition functor $\bR i_* q^*:
\db{\coh(X_0)}\to \db{\coh(Y_0)}$ and
denote it by $\Phi_{Z}$.

The aim of this paragraph is to prove the following theorem.

\begin{theorem}\label{main1}
The functor $\Phi_{Z}:\db{\coh(X_0)}\to \db{\coh(Y_0)}$
defined  by formula
$$
\Phi_{Z}(\cdot)=\bR i_* q^*(\cdot)
$$
induces a functor $\ove{\Phi}_{Z}: \dsing{X_0}\lto \dsing{Y_0}$
which is an equivalence.
\end{theorem}

\begin{remark}
{\rm An assertion  analogous to  this theorem is known in local theory of singularities
 as Kn\"orrer
periodicity. It was proved for maximal Cohen-Macalay modules
over regular
analitic k-algebra  by
Kn\"orrer (\cite{Kn}, Th.3.1).
He used a matrix factorization introduced by Eisenbud in
\cite{Eis}.}
\end{remark}

At first, we consider  partial compactifications of $Y$ and $Y_0$.
Denote by $\ove{Y}$ the scheme $X\times \PP^2$ and let
$\ove{Y}_0\subset \ove{Y}$ be a closed subscheme which is given by
equation $f z^2+xy=0$, where now $x,y,z$ are projective coordinates on
$\PP^2$. There is a flat map $\pi: \ove{Y}_0\to X$ which is
a conic bundle over $X$.
The scheme $\ove{Y}_0$ is a partial compactification
of $Y_0$ such that $\Sing(\ove{Y}_0)=\Sing(Y_0)$.
Hence, by  Proposition \ref{locpr} we have an equivalence
$\dsing{\ove{Y}_0}\cong \dsing{Y_0}$ of triangulated categories
of singularities.

Consider the cartesian square
$$
\begin{CD}
\wt{Z} @>\ove{i}>> & \ove{Y}_0\\
@VpVV     &    @VV{\pi}V\\
X_0 @>j>> &X
\end{CD}
$$
Here $\wt{Z}=X_0\times_X \ove{Y}_0$ is a fiber product.
The scheme $\wt{Z}$ is a union of two components $\ove{Z}_1\cup\ove{Z}_2$.
 Each $\ove{Z}_i$ is isomorphic
to $\PP^1\times X_0$ and their intersection is isomorphic to
$X_0$. The component $\ove{Z}_1$ is a partial compactification of
$Z=\AA^1\times X_0$. We denote by $i_1, i_2$ the closed embeddings
of $\ove{Z}_1, \ove{Z}_2$ in $\ove{Y}_0$ and we denote by
$W\subset \ove{Y}_0$ the intersection $\ove{Z}_1\cap \ove{Z}_2$.
We know that $W$ is isomorphic to $X_0$. Moreover,
$\Sing(\ove{Y}_0)$ is contained in $W$ and coincides with
$\Sing(X_0)$ under the isomorphism $W\cong X_0$.

\begin{lemma}
The closed embedding $i_1: \ove{Z}_1\hookrightarrow \ove{Y}_0$
is regular and $\ove{Z}_1$ is a Cartier divisor in $\ove{Y}_0$.
The restriction of the line bundle $\O_{\ove{Y}_0}(\ove{Z}_1)$
on $\ove{Z}_1$ is isomorphic to
$\O_{\ove{Z}_1}(-1)=\O_{\PP^1}(-1)\boxtimes \O_{X_0}$.
\end{lemma}
\begin{dok}
Consider a Cartier divisor in $\ove{Y}_0$ given by the equation $x=0$.
It is the union of $\ove{Z}_1$ and $D$, where $D$ is a component
of $\ove{Y}_0\setminus Y_0$. Hence, $D$ does not meet the singularities
of $\ove{Y}_0$. This implies that $D$ is a Cartier divisor.
Therefore, $\ove{Z}_1$ is a Cartier divisor too and $i_1$ is a regular embedding.

Consider the divisor $\wt{Z}=\ove{Z}_1\cup \ove{Z}_2$. It is equal to
$\pi^{-1}(X_0)$ and its restriction on $\ove{Z}_1$ is trivial.
The restriction of $\O(\ove{Z}_2)$ on $\ove{Z}_1$ is
$\O_{\PP^1}(1)\boxtimes \O_{X_0}$, because the intersection $\ove{Z}_2
\cap \ove{Z}_1=W$ is isomorphic to $X_0$.
Hence, the restriction of $\O(\ove{Z}_1)$ on $\ove{Z}_1$ is
$\O_{\PP^1}(-1)\boxtimes \O_{X_0}$.
\end{dok}

Denote by $p_1$ the projection of $\ove{Z}_1$ on $X_0$
and consider the commutative diagram
$$
\begin{CD}
\ove{Z}_1 @>{i_1}>> & \ove{Y}_0\\
@V{p_1}VV     &    @VV{\pi}V\\
X_0 @>j>> &X.
\end{CD}
$$
\begin{proposition}
The functor $\Phi_{\ove{Z}_1}=\bR i_{1*} p_1^*:\db{\coh(X_0)}\to \db{\coh(\ove{Y}_0)}$
 is fully faithful.
\end{proposition}
\begin{dok}
The functor $\Phi_{\ove{Z}_1}=\bR i_{1*} p_1^*$ has a right
adjoint functor $\Phi_{\ove{Z}_1 *}=\bR p_{1*} i_1^{\flat}$
where $i_1^{\flat}(\cdot)\cong \bL i_1^*(\cdot \otimes \O(Z_1))[-1]$
(see, for example, \cite{Ha} Cor.7.3).

At first, note that the functor $p_1^*:\db{\coh(X_0)}\lto
\db{\coh(\ove{Z_1})}$ is fully faithful, because by the projection
formula we have an isomorphism
$$
\bR p_{1*} p_1^*(A)\cong A \otimes \bR p_{1*} \O_{\ove{Z}_1}\cong
A\otimes \O_{X_0}\cong A
$$
for any $A\in \db{\coh(X_0)}$.

Consider the canonical transformation of functors
$\id \to \Phi_{\ove{Z}_1 *}\Phi_{\ove{Z}_1}$.
We have to show that this transformation
is an isomorphism of functors.
Take an object $A\in \db{\coh(X_0)}$ and consider an exact triangle
\begin{equation}\label{agtr}
A\lto \Phi_{\ove{Z}_1 *}\Phi_{\ove{Z}_1}(A)\lto C.
\end{equation}
To show that $C=0$ it is sufficient to check that $\bR j_* C=0$,
because $j$ is a closed embedding.

Since the functor $p_1^*$ is fully faithful, the triangle (\ref{agtr})
is the image of the triangle
\begin{equation}\label{agtr1}
p_1^* A\lto i_1^{\flat} \bR i_{1*} p_1^* A \lto B
\end{equation}
under the functor $\bR p_1^*$, where the former morphism
is a canonical map induced by the natural transformation $\id\to i_1^{\flat}
\bR i_{1*}$. Applying the functor $\bR i_{1*}$ to the exact triangle
(\ref{agtr1}) we obtain the following triangle
$$
\bR i_{1*}p_1^* A\lto \bR i_{1*} i_1^{\flat} \bR i_{1*} p_1^* A
\lto \bR i_{1*} B.
$$
This triangle is split by the canonical morphism $\bR i_{1*}
i_1^{\flat} \bR i_{1*} p_1^* A \to \bR i_{1*}p_1^* A$ which is
obtained by the tensoring product of the object $\bR i_{1*}p_1^*
A$ with a map $\alpha$ from the following exact triangle
$$
i_{1*}\O_{\ove{Z}_1}(\ove{Z}_1)[-1]\stackrel{\alpha}{\lto}
\O_{\ove{Y}_0}\lto \O_{\ove{Y}_0}(Z_1)\lto i_{1*}\O_{\ove{Z}_1}(\ove{Z}_1).
$$
Therefore, the object $\bR i_* B$ is isomorphic to
$\bR i_*p_1^* A(\ove{Z_1})=\bR i_*(p_1^* A\otimes \O_{\ove{Z}_1}(-1))$.
Now we have the following sequence of isomorphisms
\begin{multline*}
\bR j_* C\cong
\bR j_* \bR p_{1*} B\cong
\bR \pi_* \bR i_{1*} B\cong
\bR \pi_* \bR i_{1*} (p_1^* A\otimes \O_{\ove{Z}_1}(-1))\\
\cong
\bR j_* \bR p_{1*}(p_1^* A\otimes \O_{\ove{Z}_1}(-1))\cong
\bR j_* (A\otimes \bR p_{1*} \O_{\ove{Z}_1}(-1))\cong
0.
\end{multline*}
Thus the functor $\Phi_{\ove{Z}_1}=\bR i_{1*} p_1^*$ is fully
faithful.
\end{dok}

\begin{corollary}
The functor $\Phi_{\ove{Z}_1}:\db{\coh(X_0)}\to \db{\coh(\ove{Y}_0)}$
induces a functor $\ove{\Phi}_{Z}: \dsing{X_0}\lto \dsing{\ove{Y}_0}$
which is fully faithful.
\end{corollary}
\begin{dok}
The functors $p_{1}^*$ and $i_{1}^{\flat}=\bL i_1^*(\cdot\otimes
\O(\ove{Z}_1))[-1]$
take  perfect complexes to  perfect complexes as functors
of inverse images.
The functors $\bR i_{1*}$ and $\bR p_{1*}$
also preserve perfect complexes, because both
 morphisms $i_1$ and $p_1$ are finite Tor-dimension,
proper, and of finite type.

Thus, we get a functor
$\ove{\Phi}_{\ove{Z}_1}:\dsing{X_0}\to \dsing{Y_0}$ and
this functor has the right adjoint $\ove{\Phi}_{\ove{Z}_1 *}$.
As the composition
${\Phi}_{\ove{Z}_1*}{\Phi}_{\ove{Z}_1}
$ is isomorphic to the identity functor, the composition
$\ove{\Phi}_{\ove{Z}_1 *} \ove{\Phi}_{\ove{Z}_1}$ is also
isomorphic to the identity functor.
\end{dok}

Thus, we  obtained the functor
$\ove{\Phi}_{\ove{Z}_1}: \dsing{X_0}\lto \dsing{\ove{Y}_0}$
and showed that this functor is fully faithful.
To complete the proof of the theorem we have to check
that this functor is an equivalence.
We show that  any object $\E \in \dsing{\ove{Y}_0}$ satisfying
condition $\ove{\Phi}_{\ove{Z}_1 *} \E=0$ is zero object.
The property to be an equivalence will easy follow from this fact.

\begin{lemma}\label{semor}
Any object $A\in \db{\coh(\ove{Z}_1)}$ such that $\bR p_{1*}A=0$
is isomorphic to an object of the form $p_1^* B\otimes
\O_{\ove{Z}_1}(-1)$ for some $B\in \db{\coh(X_0)}$.
\end{lemma}
\begin{dok}
Denote by $B$ the object $\bR p_{1*}(A(1))$. We have a natural map
$p_1^* B\otimes \O_{\ove{Z}_1}(-1)\to A$. Denote by $C$ the cone of this map.
The object $C$ satisfies the following conditions
$$
\bR p_{1*} C=0
\qquad
\text{and}
\qquad
\bR p_{1*}(C(1))=0,
$$
the latter condition is following from the fact that $p^*_1$ is fully faithful.
Now any sheaf $\O_{\ove{Z}_1}(n)$ has a resolution of the form
$$
\O_{\ove{Z}_1}(-1)^{\oplus n}\lto \O_{\ove{Z}_1}^{\oplus (n+1)}\lto
\O_{\ove{Z}_1}(n).
$$
Tensoring this sequence with the object $C$ and applying the functor
$\bR p_{1*}$ we get that
$\bR p_{1*}(C(n))=0$ for all $n$. Hence, the object $C$ is zero, because
the sheaf $\O_{\ove{Z}}(1)$ is relatively ample.
\end{dok}

Actually, this Lemma shows us that the category
$\db{\coh(\ove{Z}_1)}$ has a semiorthogonal decomposition of the form
$$
\bigl\langle p_1^*\db{\coh(X_0)}\otimes O_{\ove{Z}_1}(-1),
\quad p^*_1
\db{\coh(X_0)}\bigl\rangle
$$
(for definition see \cite{BO}).
It can be proved for any $\PP^1$\!-bundle and moreover
for the projectivization of any bundle. The proof for smooth
base can be found in \cite{Or1} and it works for any base.

The second Lemma is also almost evident.

\begin{lemma}\label{locdiv}
Let $i: Z\hookrightarrow Y$ be a closed embedding of a Cartier
divisor. Let $\E$ be a sheaf on $Y$ such that its restriction to
the complement $U=X\setminus Z$ is locally free and $\bL i^* \E$
is isomorphic to a locally free sheaf on $Z$. Then $\E$ is locally
free on $Y$.
\end{lemma}
\begin{dok}
To prove that $\E$ is locally free it is sufficient to show that
for any closed point $t:y\hookrightarrow Y$ we have the equalities
$$
\Ext^{i}(\E, t_*\O_y)=0
$$
for all $i>0$. The sheaf $\E$ is locally free on $U$. Hence, we only need to
consider the case  $y\in Z$. This means that $t=i \cdot t'$ where $t':
y\hookrightarrow Z$ is closed embedding. In this case
$$
\Ext^{i}_{Y}(\E, t_*\O_y)=\Hom^{i}_{Z}(\bL i^* \E, t'_*\O_y)=0
$$
for $i>0$, because $\bL i^* \E$ is isomorphic to a locally free sheaf on $Z$.
\end{dok}

Using these two lemmas we can prove the following proposition.

\begin{proposition}\label{zerob}
Assume that an object $\E\in \dsing{\ove{Y}_0}$ satisfies the condition
$\ove{\Phi}_{\ove{Z}_1 *} \E=0$. Then $\E=0$ in $\dsing{\ove{Y}_0}$.
\end{proposition}
\begin{dok}
At first, note that all schemes $X_0, \ove{Z}_1, \ove{Y}_0$
are Gorenstein. By Proposition \ref{fint} we can assume that
$\E$ is a sheaf and $\ext^{i}(\E, \O_X)=0$ for all $i\ne 0$.
Note that any such $\E$ is locally free on the complement
$X\setminus \Sing(X)$.
In addition, for such $\E$ we have $\bL i_1^* \E\cong i_1^*\E$ is a sheaf.

Denote by $\L$ the relatively ample line bundle on $\ove{Y}_0$
obtained by the restriction of the line bundle
$\O_{\ove{Y}}(1)$ on $\ove{Y}=\PP^2\times X$.
We can see that
the object $\E\otimes \L^{\otimes n}$ is isomorphic to $\E$ in the category
$\dsing{\ove{Y}_0}$,
because there is an inclusion $\E\to\E\otimes \L^{\otimes n}$
such that the support of the cokernel does not intersect
$\Sing(\ove{Y}_0)$ and, hence, it is perfect by Lemma \ref{persup}.

Take the sheaf $i_1^* \E$ and denote by $\F$
the sheaf $p_{1*}i_1^* \E$ on $X_0$.
Tensoring $\E$ with $\L^{n}$
if it is  necessary  we can reduce the situation to the case
$$
\bR^i p_{1*} i_1^* \E=0
\quad\text{for all}\quad
i>0
\quad
\text{and}
\quad
p_1^* \F\stackrel{\alpha}{\lto} i^*_1\E
\quad
\text{is surjective.}
$$
By Lemma \ref{semor} the kernel of $\alpha$ is a sheaf of the form
$p_1^*\G\otimes\O_{\ove{Z}_1}(-1),$ where $\G$ is a sheaf on
$X_0$. By the assumption, the sheaf $\F$ is  perfect as a complex
on $X_0$. Moreover, the sheaf $\G\cong p_{1*}i_1^*(\E\otimes
\L^{-1})$ is also perfect as a complex on $X_0$. Hence, the sheaf
$i^*_1\E$ on $\ove{Z}_1$ is a perfect complex on $\ove{Z}_1$. On
the other hand, we know that the sheaf $\E$ has a right locally
free resolution. Hence the sheaf $i_1^* \E$ as the restriction
with respect to a regular embedding of a divisor also has a right
locally free resolution. By Lemma \ref{rres} this implies that
$$
\ext^i(i^*_1 \E, \O_{\ove{Z}_1})=0
\quad
\text{for all}
\quad
i>0.
$$
By Lemma \ref{locfr} this means that $i^*_1 \E$ is locally
free on $\ove{Z}_1$. Now apply Lemma \ref{locdiv} we get that
$\E$ is locally free on whole $\ove{Y}_0$.
Hence, it is isomorphic to the zero object in $\dsing{\ove{Y}_0}$.
\end{dok}

\noindent{\bf Proof of Theorem \ref{main1}}
We already know that $\ove{\Phi}_{\ove{Z}_1}$ is fully faithful.
Take an object $A\in \dsing{\ove{Y}_0}$ and consider the natural map
$ \ove{\Phi}_{\ove{Z}_1 }\ove{\Phi}_{\ove{Z}_1 *} A\to A$.
Denote
by $C$ its cone. Apply the functor $\ove{\Phi}_{\ove{Z}_1 *}$
to the obtained exact triangle.
Since $\ove{\Phi}_{\ove{Z}_1}$ is fully faithful,
we get that $\ove{\Phi}_{\ove{Z}_1 *}C=0$. By Proposition
\ref{zerob} the object $C$ is the zero object. Hence, the functor
$\ove{\Phi}_{\ove{Z}_1 *}$ is also fully faithful and, consequently,
it is an equivalence.
It remains only to note  that the functor
$\ove{\Phi}_{Z}: \dsing{X_0}\to \dsing{Y_0}$ is the composition
of the functor $\ove{\Phi}_{\ove{Z}_1}:\dsing{X_0}\to
\dsing{\ove{Y}_0}$ and the functor
$\ove{J}^*:\dsing{\ove{Y}_0}\to \dsing{Y_0}$, where
$J:Y_0\hookrightarrow \ove{Y}_0$ is the open embedding.
Both of these functors are equivalence. Hence, $\ove{\Phi}_{Z}$
is an equivalence too.
\hfill$\Box$

\section{Triangulated category of D-branes of type B in Landau-Ginzburg model}
\label{three}

\subsection{Kontsevich's proposal and categories of pairs}

A mathematical definition of the categories of D-branes of type B in
Landau-Ginzburg models is proposed by M.Kontsevich (see also
\cite{KL}).

By a Landau-Ginzburg model we mean the following data: a smooth
variety (or regular scheme) $X$ and a regular function $W$ on $X$
such that the morphism $W: X\lto \AA^1$ is flat (for the
definition of B-branes we don't need a symplectic form on $X$
which have to be in a LG model too).

With any point $w_0\in \AA^1$ we can associate a differential
$\ZZ/2\ZZ$\!-graded category $\dg_{w_0}(W)$, an exact category
$\abe_{w_0}(W)$, and a triangulated category $\dt_{w_0}(W)$. We
give constructions of these categories under the condition that
$X=\Spec(A)$ is affine (see \cite{KL}). The general definition is
more sophisticated.

Since the category of coherent sheaves on an affine scheme
$X=\Spec(A)$ is the same as the category of finitely generated
$A$\!-modules we will frequently go from sheaves to modules and
back. Note that under this equivalence locally free sheaves are
the same as projective modules.

Objects of all these categories are  ordered pairs
$$
\ove{P}:=\Bigl(
\xymatrix{
P_1 \ar@<0.6ex>[r]^{p_1} &P_0 \ar@<0.6ex>[l]^{p_0}
}
\Bigl)
$$
where $P_0, P_1$ are  finitely generated
projective $A$\!-modules and
the compositions
${p_0 p_1}$ and ${p_1 p_0}$
are the multiplications by the element $(W-w_0)\in A$.

Morphisms from $\ove{P}$ to $\ove{Q}$
in the category $\dg_{w_0}(W)$ form $\ZZ/2\ZZ$\!-graded
complex
$$
{\Hhom}(\ove{P},\; \ove{Q})=\bigoplus_{i,j} \Hom(P_i, Q_j)
$$
with a natural grading $(i-j)\mod 2,$ and with a
differential $D$ acting on  homogeneous elements
of degree $k$ as
$$
D f=q \comp f-(-1)^{k}f\comp p .
$$

The space of morphisms $\Hom(\ove{P},\; \ove{Q})$
in the category $\abe_{w_0}(W)$
is the space  of homogeneous of degree 0
morphisms in $\dg_{w_0}(W)$ which commutes with the differential.
The space of morphisms in  the category $\dt_{w_0}(W)$ is
the space of morphisms in $\abe_{w_0}(W)$
 modulo null-homotopic morphisms, i.e.
$$
\Hom_{\abe_{w_0}(W)}(\ove{P},\; \ove{Q})= \rZ^0(\Hhom(\ove{P},\; \ove{Q})),
\qquad
\Hom_{\dt_{w_0}(W)}(\ove{P},\; \ove{Q})= \rH^0(\Hhom(\ove{P},\; \ove{Q})).
$$
Thus a morphism $f:\ove{P}\to\ove{Q}$ in the category
$\abe_{w_0}(W)$ is a pair of morphisms
$f_1: P_1\to Q_1$ and $f_0: P_0\to Q_0$ such that
$f_1p_0=q_0 f_0$ and $q_1f_1=f_0 p_1$.
The morphism $f$ is null-homotopic if there are two morphisms
$s: P_0\to Q_1$ and $t:P_1\to Q_0$ such that
$f_1=q_0 t + s p_1$ and $f_0=t p_0 + q_1 s$.

It is clear that the category $\abe_{w_0}(W)$
is an exact category with respect to componentwise
monomorphisms and epimorphisms (see definition in \cite{Qu}).

\begin{remark}{\rm
The remarkable fact is that such construction appeared many years
ago in the paper \cite{Eis} and is known for specialist in
singular theory as a matrix factorization.}
\end{remark}

The category $\dt_{w_0}(W)$ can be
endowed with a natural structure of a triangulated category.
To determine it we have to define a translation functor $[1]$
and a class of exact triangles.

The translation functor can be defined as a functor
that takes an object $\ove{P}$ to the object
$$
\ove{P}[1]=\xymatrix{
\Bigl(
P_0 \ar@<0.6ex>[r]^{-p_0} &P_1 \ar@<0.6ex>[l]^{-p_1}
}
\Bigl),
$$
i.e. it changes the order of the modules and signs of the morphisms,
and takes a morphism $f=(f_0, f_1)$ to the morphism
$f[1]=(f_1, f_0)$.
We see that the functor $[2]$ is the identity functor.

For any morphism $f:\ove{P}\to\ove{Q}$ from the category
$\abe_{w_0}(W)$ we define a mapping cone $C(f)$ as an object
$$
C(f)=\xymatrix{
\Bigl(
Q_1\oplus P_0 \ar@<0.6ex>[r]^{c_1} & Q_0\oplus P_1 \ar@<0.6ex>[l]^{c_0}
}
\Bigl)
$$
such that
$$
c_0=
\begin{pmatrix}
q_0 & f_1\\
0 & -p_1
\end{pmatrix},
\qquad
c_1=
\begin{pmatrix}
q_1 & f_0\\
0 & -p_0
\end{pmatrix}.
$$
There are  maps
$g: \ove{Q}\to C(f), \; g=(\id , 0)$ and $h: C(f)\to\ove{P}[1],\;
h=(0, -\id)$.

Now we define a standard triangle in the category
$\dt_{w_0}(W)$ as a triangle of the form
$$
\ove{P}\stackrel{f}{\lto} \ove{Q}\stackrel{g}{\lto} C(f)
\stackrel{h}{\lto} \ove{P}[1].
$$
for some $f\in \abe_{w_0}(W)$.

\begin{definition}
A triangle
$\ove{P}{\to} \ove{Q}{\to} \ove{R}
{\to} \ove{P}[1]$ in $\dt_{w_0}(W)$
will be  called  an exact triangle if it is isomorphic to a
standard triangle.
\end{definition}

\begin{proposition}\label{trstr}
The category $\dt_{w_0}(W)$ endowed with the translation functor $[1]$
and the above class of exact triangles becomes a triangulated category.
\end{proposition}
\begin{dok}The proof is the same as the proof of the analogous
result for a usual homotopic category (see, for example, \cite{GM}
or \cite{KS}).
\end{dok}

\begin{definition} We define a category of D-bran of type B (B-branes) on $X$ with the
superpotential $W$ as the product $\dt(W)=\mathop{\prod}\limits_{w\in \AA^1}
\dt_{w}(W)$.
\end{definition}
 Note that since $X$
is regular, the set of points on $\AA^1$ with singular fibers  is finite
\cite{H}, III,Cor.10.7 (we suppose here that we work over a field of
characteristic 0). It will be shown that the category $\dt_{w}(W)$ is trivial
if the fiber over point $w$ is smooth. Hence, the category $\dt(W)$ is a
product of finitely many numbers of categories.

The product of two triangulated categories $\D_1$ and $\D_2$ is the same as
their orthogonal sum. Objects of it are pairs $(A, B)$, where $A\in\D_1$ and
$B\in\D_2$. The space of morphisms from $(A, B)$ to $(A', B')$ is the sum
$\Hom(A,A')\oplus\Hom(B, B')$. The translation functor and exact triangles are
defined by componentwise.

The construction of the category $\dt_{w_0}(W)$ fits into a
general construction of the {\sf stable category} associated with
an exact category. We briefly recall the definition of an exact
category, which was introduced by Quillen in \cite{Qu}.

An exact category $\E$ is an additive category together with a
choice of a class of sequences $\{ F\mon E \epi G \}$ said to be
exact. This determines two classes of morphisms: the admissible
epimorphisms $E\epi G$ and the admissible monomorphisms $F\mon E$.
The exact category is to satisfy the following axioms: Any
sequence isomorphic to an exact sequence is exact. In any exact
sequence $F\stackrel{i}{\mon} E\stackrel{p}{\epi} G$, the map $i$
is a kernel of $p$ and $p$ is a cokernel of $i$. The class of
admissible monomorphisms is closed under composition and is closed
under cobase change by pushout along an arbitrary map $F\to F'$.
Dually, the class of admissible epimorphisms is closed under
composition and under base change by pullback along an arbitrary
map $G'\to G$.
This definition is equivalent to the original definition of Quillen
(see \cite{Ke1}).

An object $I\in \E$ is {\sf injective} (resp. $P$ is {\sf projective})
if the sequence
$$
\Hom(E, I)\to \Hom(F, I)\to 0
\qquad
(\text{resp.}
\;
\Hom(P, E)\to \Hom(P, G)\to 0)
$$
is exact for each admissible monomorphism (resp. epimorphism).
We say that $\E$ has enough injectives, if for each $F\in \E$
there is an admissible monomorphism in an injective.

If $\E$ has enough injectives then we can define the
{\sf stable category} $\underline{\E}$ as a category
which has the same objects as $\E$ and a morphism in $\underline{\E}$ is the
equivalence class $\ove{f}$ of a morphism $f$ of $\E$
modulo the subgroup of morphisms factoring through an injective in $\E$
(see, for example, \cite{Hap, Ke, Ke1}).

If  $\E$ also has enough projectives
(i.e. for each $G\in\E$ there is
an admissible epimorphism $P\epi G$ with projective $P$),
and
the classes of projectives and injectives coincide, then
$\E$ is called a {\sf Frobenius category}.
The stable category
$\underline{\E}$ associated with a Frobenius category
has a natural structure of a triangulated category (\cite{Hap}).

In our case, the category $\abe_{w_0}(W)$ is an exact category
with respect to  componentwise monomorphisms and epimorphisms.
Moreover, it can be shown that $\abe_{w_0}(W)$ is a Frobenius
category and the class of injectives consists exactly of
 homotopic to zero pairs.
Hence, the category $\dt_{w_0}(W)$ is nothing more than the stable
category associated to the exact category $\abe_{w_0}(W)$
and has the natural structure of a triangulated category.
It gives another proof of Proposition \ref{trstr}.

\subsection{Categories of pairs and categories of singularities}

In this paragraph we establish a connection between category of
B-branes in a Landau-Ginsburg model and triangulated categories of
singular fibres, introduced above.
As it turned out this
connection for maximal Cohen-Macaulay modules over local ring appeared
 in the paper \cite{Eis}, Sect.6.

Further we suppose $w_0=0$.
Denote by $X_{0}$ the fiber of $f: X\to \AA^{1}$ over
the point $0$.
With any pair $\ove{P}$ we can associate a short exact sequence
\begin{equation}\label{shseq}
0\lto P_1 \stackrel{p_1}{\lto} P_0 \lto \Coker p_1 \lto 0.
\end{equation}

We can attach to an object $\ove{P}$ the sheaf $\Coker p_1$.
This is a sheaf on $X$. But the multiplication by $W$
annihilates it. Hence, we can consider $\Coker p_1$
as a sheaf on $X_{0}$. Any morphism $f:\ove{P}\to\ove{Q}$
in $\abe_{0}(W)$ gives a morphism between
 cokernels.
This way we get a functor $\Cok:\abe_{0}(W)\lto \coh(X_{0})$.

\begin{lemma}
The functor $\Cok$ is full.
\end{lemma}
\begin{dok}
Any map $g: \Coker p_1 \to \Coker q_1$ can be extended to a map
of exact sequences
$$
\begin{CD}
0 @>>> P_1 @>p_1>> P_0 @>>>  \Coker p_1 @>>> 0\\
&& @V{f_1}VV @VV{f_0}V @VV{g}V \\
0 @>>> Q_1 @>q_1>> Q_0 @>>>  \Coker q_1 @>>> 0,
\end{CD}
$$
because $P_1$ and $P_0$ are projective.
To prove the lemma it is sufficient  to show that $f=(f_1, f_0)$ is a map of pairs,
i.e $f_1p_0=q_0f_0$.
We have the sequence of equalities
$$
q_1(f_1p_0-q_0f_0)=f_0 p_1 p_0 - q_1q_0 f_0=
f_0 W- W f_0=0.
$$
Since $q_1$ is an embedding, we obtain that $f_1p_0=q_0f_0$.
\end{dok}

\begin{lemma}\label{extv}
For any pair $\ove{P}$ the coherent sheaf
$\Coker p_1$ on $X_0$ satisfies the condition
$$
\Ext^{i}(\Coker p_1, \O_{X_0})=0
$$
for all $i>0$.
\end{lemma}
\begin{dok}
First note that $X_{0}$ is Gorenstein as a full intersection.
Consider the restriction of the sequence (\ref{shseq}) on $X_0$.
We obtain an exact sequence
$$
0\to \Coker p_1\lto P_1/W \stackrel{p_1|_{W}}{\lto} P_0/W \lto \Coker p_1 \to 0.
$$
This gives us a periodic resolution of the sheaf $\Coker p_1$
$$
0\to \Coker p_1\lto P_1/W \stackrel{p_1|_{W}}{\lto} P_0/W \lto
P_1/W\lto\cdots.
$$
Since $X_0$ is affine and Gorenstein the existence of such
resolution implies that
$$
\Ext^{i}(\Coker p_1, \O_{X_0})=0
$$
for all $i>0$ ( by Lemma \ref{rres}).
\end{dok}

Now we show that the functor $\Cok$ induces an exact functor between
triangulated categories $\dt_{0}(W)$ and $\dsing{X_{0}}$.

\begin{proposition} There is a functor $F$ which completes
the following commutative diagram
$$
\begin{CD}
\abe_{0}(W) @>\Cok>> & \coh(X_{0})\\
@VVV&@VVV\\
\dt_{0}(W) @>F>>& \dsing{X_{0}}.
\end{CD}
$$
Moreover, the functor $F$ is an exact functor
between triangulated categories.
\end{proposition}
\begin{dok}
We have a functor  $\abe_{0}(W)\to \dsing{X_0}$ which is the
composition of $\Cok$ and the natural functor from $\coh(X)$ to
$\dsing{X_0}$. To prove the existence of a functor $F$ we need to
show that any morphism $\ove{f}=(f_1, f_0): \ove{P} \to \ove{Q}$
which is homotopic to $0$ goes to $0$\!-morphism in $\dsing{X_0}$.
Fix a homotopy $(t, s)$ where $t:P_1\to Q_0$ and $s:P_0\to Q_1$.
Consider the following decomposition of $\ove{f}$:
$$
\xymatrix{
P_1 \ar[d]_{(t, f_1)} \ar@<0.6ex>[r]^{p_1}  & P_0
\ar@<0.6ex>[l]^{p_0} \ar[d]^{(s, f_0)}\ar[r] & \Coker p_1 \ar[d]
\\
Q_0\oplus Q_1 \ar[d]_{pr} \ar@<0.6ex>[r]^{c_1} & Q_1\oplus
Q_0\ar[d]^{pr} \ar@<0.6ex>[l]^{c_0} \ar[r] & Q_0/W \ar[d]
 & {}\save[]*\txt{where}\restore&
{c_0=
\begin{pmatrix}
-q_0 & \id\\
0 & q_1
\end{pmatrix},
\;
c_1=
\begin{pmatrix}
-q_1 & \id\\
0 & q_0
\end{pmatrix}.
}
\\
Q_1 \ar@<0.6ex>[r]^{q_1} & Q_0 \ar@<0.6ex>[l]^{q_0}\ar[r] &
\Coker q_1
}
$$
This gives the decomposition of $F(\ove{f})$ through
a locally free object $Q_0/W$.
Hence, $F(\ove{f})=0$ in the category $\dsing{X_0}$.
We leave to the reader the proof that the functor $F$ is exact.
\end{dok}

\begin{lemma}\label{embob}
If $F\ove{P}=0$ then $\ove{P}=0$ in $\dt_{0}(W)$.
\end{lemma}
\begin{dok}
If $F\ove{P}=0$ then $\Coker p_1$ is a perfect complex.
This implies that it is locally free by Lemmas \ref{locfr} and \ref{extv}.
Hence, there is a map $f: \Coker p_1 \to P_0/W$ which splits
the epimorphism $pr: P_0/W\to\Coker p_1$.
It can be lifted to a map from $\{P_1\stackrel{p_1}{\lto} P_0 \}$
to $\{P_0\stackrel{W}{\lto}P_0\}$.
Consider a diagram
$$
\begin{CD}
P_1 @>p_1>> & P_0 @>>> & \Coker p_1\\
@V{t}VV & @VV{u}V & @VV{f}V\\
P_0 @>W>> & P_0 @>>> &P_0/W\\
@V{p_0}VV & @VV{\id}V & @VV{pr}V\\
P_1 @>p_1>> & P_0 @>>> & \Coker p_1.
\end{CD}
$$
Since the composition $pr f$ is identical, the map $(p_0 t, u)$
from the pair $\{ P_1\stackrel{p_1}{\lto} P_0 \}$ to itself is
homotopic to the identical map of this pair. Hence, there is a map
$s: P_0\to P_1$ such that
$$
\id_{P_1} - p_0 t= s p_1 \qquad \text{and} \qquad p_1
s=\id_{P_0}-u.
$$
Moreover, we have the following equalities
$$
0= (u p_1- W t)=(u p_1- tW)=(u -t p_0) p_1.
$$
This gives us that $(u-t p_0)=0$, because no maps from $\Coker
p_1$ to $P_0$. Finally, we obtain two morphisms $t$ and $s$ such
that
$$
\id_{P_1}= p_0 t + s p_1 \qquad \text{and} \qquad \id_{P_0}=p_1 s
+ t p_0.
$$
Hence, the pair $\ove{P}$ is isomorphic to the zero object in the
category $\dt_{0}(W)$.
\end{dok}

\begin{theorem}\label{main2}
The functor $F:\dt_{0}(W)\lto \dsing{X_{0}}$ is an exact
equivalence.
\end{theorem}
\begin{dok}
By Lemma \ref{extv} the coherent sheaves $\Coker p_1$ and $\Coker q_1$
satisfy the condition of  Proposition \ref{stab} with $N=0$.
This gives an isomorphism
$$
\Hom_{\dsing{X_{0}}} (\Coker p_1 , \Coker q_1)\cong
\Hom_{\coh(X_{0})}  (\Coker p_1 , \Coker q_1)/\R
$$
where $\R$ is the subspace of morphisms factoring through a
locally free sheaf.
Since $\Cok$ is full we get that the functor $F$ is also full.

Now we prove that $F$ is faithful. It is a standard consideration.
Let $f: \ove{P}\to \ove{Q}$ be a morphism for which $F(f)=0$.
Suppose that $f$ sits in an exact triangle
$$
\ove{P}\stackrel{f}{\lto}\ove{Q}\stackrel{g}{\lto}\ove{R}
$$
Then the identity map of $F\ove{Q}$ factors through the map
$F\ove{Q}\stackrel{Fg}{\lto} F\ove{R}$.
Since $F$ is full, there is a map $h:\ove{Q}\to \ove{Q}$
factoring through $g: \ove{Q}\to \ove{R}$ such that $Fh=\id$.
Hence, the cone $C(h)$ of map $h$ goes to zero under the functor $F$.
By Lemma \ref{embob} the object $C(h)$ is the zero object as well,
so $h$ is an isomorphism.
Thus $g:\ove{Q}\to\ove{R}$ is a split monomorphism and $f=0$.

To complete the proof that $F$ is an equivalence we need to show that
every object $A\in \dsing{X_0}$ is isomorphic to $F\ove{P}$ for some $\ove{P}$.
By Proposition \ref{fint} any object $A\in \dsing{X_0}$ is isomorphic to
the image of a coherent sheaf $\F$ such that $\ext^{i}(\F, \O_X)=0$
for all $i>0$.
Consider an epimorphism $P_0\to \F$ of sheaves on $X$  with
locally free $P_0$. Denote by $p_1: P_1\to P_0$ the kernel of this map.
 Since the multiplication on $W$ gives the zero map on $\F$,
there is a map $p_0: P_0\to P_1$ such that
$p_0 p_1=W$ and $p_1 p_0=W$. We get a pair
$$
\ove{P}:=\Bigl(
\xymatrix{
P_1 \ar@<0.6ex>[r]^{p_1} &P_0 \ar@<0.6ex>[l]^{p_0}
}
\Bigl)
$$
And we need only to check that $P_1$ is locally free. It  follows
from the fact that for any closed point $t:x\hookrightarrow
X$ we have
\begin{equation}\label{last}
\Ext^i(P_1, t_* \O_x)= 0
\end{equation}
for all $i>0$. To show it we note that by Lemma \ref{rres} the
sheaf $\F$ has a right locally free resolution on $X_0$.
For any local free sheaf $Q$ on $X_0$ we have
$
\Ext^i_X(Q, t_* \O_x)= 0
$
for $i>1$. And since the category of coherent sheaves on $X$
has finite cohomological dimension we obtain    $\Ext^i_X(\F, t_* \O_x)= 0$
for $i>1$. This immediately implies (\ref{last}).
\end{dok}

\begin{corollary} The category of B-branes on
smooth $X$ with a superpotential $W$ is equivalent to the product
$\mathop{\prod}\limits_{w\in \AA^1} \dsing{X_{w}}$, and this
product is finite.
\end{corollary}

\subsection{Some simple calculations}
In this paragraph we give a description of the category of
B-branes in the Landau-Ginzburg model with a superpotential
$W'=z_0^n + z_1^2+\cdots + z_{2k}^2$ which is given on
$\CC^{2k+1}$. (The descriptions of these categories and another
categories which is connected with another Dynkin diagram are
known and can be found in the papers \cite{A, AR, DW}, where the
technique of Auslander-Reiten sequences is used.)

 The superpotential $W$ has only one singular point over $0$. By
Theorem \ref{main2} the category of B-branes $\dt_0(W')$ is
equivalent to the triangulated category of singularities
$\dsing{Y_0}$, where $Y_0$ is the fiber over $0$, i.e $Y_0$ is
given by the equation $W'=0$. By Theorem \ref{main1} this category
is equivalent to the category $\dsing{X_0}$, where
$X_0=\Spec(\CC[z]/z^n)$ is the fiber over $0$ of the
superpotential $W=z^n$ on $\CC$. Thus, it is enough to describe
the category $\dsing{X_0}$.

{\sc Objects.} Denote by $A$ the algebra $\CC[z]/z^n$. By
Proposition \ref{fint} any object is coming from a
finite-dimensional module $M$ over the algebra $A$. Each
$A$\!-module is nothing more as a vector space with an operator
$L$ such that $L^n=0$. Hence any module is a direct sum of the
modules $V_i$ for $i=1,...,n$, where $V_i=A/z^i$. This is the
Jordan decomposition in Jordan blocks. Moreover, the module
$V_n=A$ is  free, hence it is equal to the zero object in
$\dsing{X_0}$. So indecomposable objects in the category
$\dsing{X_0}$ are
$$
V_1, V_2, ...., V_{n-1}.
$$
All other objects are finite direct sums of $V_{\mu},\; \mu=1,... n-1$.

{\sc Morphisms.} For each pair  $V_{\mu}, V_{\nu}$
we fix a morphism
$$
{}_{\nu}\alpha_{\mu}: V_{\mu}\to V_{\nu}
$$
which is coming from the natural projection if $\mu\ge\nu$
and from the injection that sends $1\in V_{\mu}=A/z^{\mu}$ to
$z^{\nu-\mu}$ if $\nu\ge \mu$.
All other morphisms are  linear combination of compositions
of ${}_{\nu}\alpha_{\mu}$.
There are the following relations:
\begin{equation}\label{rela}
\begin{array}{lll}
1)&{}_{\mu}\alpha_{\mu}=id_{\mu},&\\
2)& {}_{\nu}\alpha_{\lambda}{}_{\lambda}\alpha_{\mu}={}_{\nu}\alpha_{\mu}
&\text{if}
\quad
\nu\ge\lambda\ge\mu
\quad
\text{or}
\quad
\nu\le\lambda\le\mu,\\
3)&{}_{\nu}\alpha_{\lambda}{}_{\lambda}\alpha_{\mu}=0
&\text{if}
\quad
\lambda\ge\mu+\nu
\quad
\text{or}
\quad
\lambda+n\le \mu+\nu,\\
4)&
{}_{\nu}\alpha_{\lambda}{}_{\lambda}\alpha_{\mu}=
{}_{\nu}\alpha_{\kappa}{}_{\kappa}\alpha_{\mu}
&
\text{if}
\quad
\lambda+\kappa=\mu+\nu.
\end{array}
\end{equation}
Using this relation we can see that the space
$\Hom(V_{\mu}, V_{\nu})$ has a basis formed by the morphisms of the form
$_{\nu}\alpha_{\lambda}{}_{\lambda}\alpha_{\mu}$, where
$\max(\mu, \nu)\le \lambda<\mu+\nu$.
Denote by $\depth V_{\mu}$ the integer number equals to $\min(\mu, n-\mu)$.
We obtain  that
$$
\dim \Hom(V_{\mu}, V_{\nu})=\min(\depth V_{\mu}, \depth V_{\nu}).
$$
Moreover, the ring $\End (V_{\mu})$ is isomorphic to
$\CC[x]/x^d$, where $d=\depth V_{\mu}$.

{\sc Translation functor.} To see the translation functor it is
convenient to draw the following pictures
$$
\xymatrix{
\buildrel{V_1}\over{\circ}&\buildrel{V_2}\over{\circ}&\cdots&
\buildrel{V_k}\over{\circ}
&\buildrel{V_1}\over{\circ}&\cdots&\buildrel{V_{k-1}}\over{\circ}&\\
-\ar@{--}[rrr]&&&&-\ar@{--}[rrr]&&& \buildrel{V_k}\over{\circ}\\
\buildrel{V_{n-1}}\over{\circ}&\buildrel{V_{n-2}}\over{\circ}&\cdots
&\buildrel{V_{k+1}}\over{\circ}&
\buildrel{V_{n-1}}\over{\circ}&\cdots&\buildrel{V_{k+1}}\over{\circ}&\\
&&{}_{n=2k+1}&&&{}_{n=2k}
}
$$

The translation functor $[1]$  is the reflection with respect to
the horizontal line, i.e. it sends $V_{\mu}$ to $V_{n-\mu}$ and
takes ${}_{\nu}\alpha_{\mu}$ to ${}_{(n-\nu)}\alpha_{(n-\mu)}$. It
is easy to see that all relations (\ref{rela}) are preserved.

{\sc Exact triangles.}
For the convenience, we will write $V_{-\mu}$ instead $V_{n-\mu}$
considering all integers modulo $n$. We also put $V_0=0$.
At first, any triangle of the form
\begin{equation}\label{fst}
V_{\mu}\stackrel{{}_{\nu}\alpha_{\mu}}{\llto}
V_{\nu}\stackrel{{}_{(\nu-\mu)}\alpha_{\nu}}{\llto}
V_{(\nu-\mu)}\stackrel{h}{\llto}
V_{(-\mu)},
\end{equation}
where $h={}_{(-\mu)}\alpha_{(\nu-\mu)}$ if $\nu-\mu\ge 0$
and $h=-{}_{(-\mu)}\alpha_{(\nu-\mu)}$ if $\nu-\mu<0$,  is exact.
If now $f: V_{\mu}\to V_{\nu}$ is  another morphism
which is a composition of some $\alpha$\!'s then using the relations
(\ref{rela}) we can represent it as ${}_{\nu}
\alpha_{\lambda}{}_{\lambda}\alpha_{\mu}$
where $\lambda>\max(\mu,\nu)$ and $\lambda<\mu+\nu$.
In this case the triangle of the form
\begin{equation}\label{lst}
V_{\mu}\stackrel{f}{\lto}
V_{\nu}\stackrel{g}{\lto}
V_{(\lambda-\mu)}\oplus V_{(\nu-\lambda)}\stackrel{h}{\lto}
V_{(-\mu)}
\end{equation}
is exact, where $g=({}_{(\lambda-\mu)}\alpha_{\nu}
\;,{}_{(\nu-\lambda)}\alpha_{\nu})^t$
and $h=
({}_{(-\mu)}\alpha_{(\lambda-\mu)},
-{}_{(-\mu)}\alpha_{(\nu-\lambda)})$.
Note that the triangle (\ref{fst}) is a particular case
of the triangle (\ref{lst}) for $\lambda=\max(\mu, \nu).$

All other exact triangles are isomorphic to the triangles defined above.

\section*{Acknowledgments}
I am grateful to Alexei Bondal, Anton Kapustin, Ludmil Katzarkov,
Bernard Keller, Maxim Kontsevich, and Leonid Posicelsky for very
useful discussions. I also thank Nikolai Tyurin for reading of a
preliminary draft of the paper and making a number of valuable
comments. I am grateful to Duco van Straten who drew my attention
to papers \cite{Eis, Kn, Bu}, he pointed me that some results of
the paper are closely related to
known results in local theory of singularities.
 For instance, the assertion of Theorem \ref{main1} is known as
Kn\"orrer periodicity and is proved for maximal Cohen-Macalay
modules over local ring in \cite{Kn}, the construction of the
section \ref{three} is also known in the local theory of
singularities as a matrix factorization, which is due to Eisenbud
\cite{Eis}, where using this construction MCM modules over local
rings were described. I also would like to express my gratitude to
Amnon Neeman who pointed out some inexactitudes in the previous
version of the paper. I would like to thank Max-Planck-Institut
f\"ur Mathematik for the hospitality during the writing of this
paper.

Finally, I want to thank my non-formal adviser Andrei Nikolaevich
Tyurin without whom this work could not be written at all.

\end{document}